\input amstex
\documentstyle{amsppt}
\loadbold
\topmatter
\title The Fine Structure of the Kasparov Groups I: \\
 Continuity of the KK-Pairing   \\
      \endtitle
 \rightheadtext{Fine Structure I }
\author    Claude L. Schochet
  \endauthor
 \affil
 Mathematics Department   \\
 Wayne State University    \\
Detroit, MI 48202      \\ \\
 \endaffil
\address       Wayne State University     \endaddress
\email     {    claude\@math.wayne.edu }     \endemail
\date {revised, January, 2001   } \enddate
\dedicatory {to Calvin C. Moore} \enddedicatory
\endkeywords
  \endsubjclass
\abstract {
 In this paper it is demonstrated
 that the Kasparov  pairing is   continuous with
respect to the natural   topology
on the Kasparov groups,
  so that a $KK$-equivalence
is an isomorphism of topological groups.
In addition, we demonstrate that the groups have a natural pseudo\-polonais
structure, and we prove that various $KK$-structural maps are continuous.

}\endabstract

\endtopmatter

\document
\magnification = 1200
\pageno=1

\def\tensor{\otimes}

\def\KK #1.#2 { KK_*( #1,#2) }
 \def\KKgraded #1.#2.#3 { KK_{#1}( #2, #3) }
 \def\bKKgraded #1.#2.#3 {\overline{ KK}_{#1}( #2, #3) }

  \def\invlim{\underset \longleftarrow\to {lim} \, }

\def\invlimone{\underset \longleftarrow\to {lim^1}  \, }

\def\ext #1.#2 {Ext _{\Bbb Z}^1( {#1} , {#2} ) }
\def\pext #1.#2 {Pext _{\Bbb Z}^1( {#1} , {#2} ) }
\def\hom #1.#2 {Hom_{\Bbb Z} (#1 ,  #2  ) }

\newpage\beginsection {1. Introduction    }

This is the first of a series of papers in which the topological
structure of the Kasparov $KK$-groups
  is developed and put to use. The Kasparov
groups $\KKgraded *.A.B $, defined for separable $C^*$-algebras $A$ and
$B$, have been shown to be powerful tools in the analysis of a
wide variety of  problems in functional analysis and in topology.
It is our hope that the \lq\lq fine structure\rq\rq\,\,    of these groups
will be of additional help and, particularly, that it will be
of use in the classification of separable nuclear simple   $C^*$-algebras.

The topological structure of the Kasparov groups was
 first studied in depth by
N. Salinas \cite{19}. We shall demonstrate that
 $\KKgraded *.{-}.{-} $ is a bifunctor to graded
pseudo\-polonais groups. The key result in this paper, Theorem
6.8, asserts
 that the $KK$-product is  jointly
continuous with respect to this topology, provided that 
the $C^*$-algebras that appear in the first variable are
$K$-nuclear. This theorem implies that a $KK$-equivalence is a
homeomorphism. These theorems are applied in \cite {25}
  to the
study of relative quasidiagonality. Central to that study
is the fine structure subgroup of $\KKgraded *.A.B $,
namely the group
$$
\pext {K_*(A)}.{K_*(B)} \,\,\cong\,\, \invlimone\hom{K_*(A_i)}.{K_*(B)}
$$
which under bootstrap hypotheses is the closure of $0$ in
$\KKgraded *.A.B$ \cite {24, 22, 23}.

The paper is organized as follows. In \S 2 we review the topology
of function spaces of $*$-homomorphisms and their quotients modulo
homotopy. Along the way we introduce K. Thomsen's useful notion of
{\it quasi-unital} maps, and we show how to use these maps to
better understand the functorial nature of $KK$. Section 3 is
devoted to showing that the $KK$-pairing is continuous with
respect to the topology on the $KK$-groups inherited from the
Cuntz representation of $KK$ as $\KKgraded 0.A.B \cong [qA,
B\tensor\Cal K]$.   In \S 4 we introduce the topology used by
Salinas as well as a natural topology associated to the Zekri
picture of $KK_1$ and we show that for $A$ $K$-nuclear these all
agree. As a consequence we show that
the $KK$-pairing is separately continuous with respect to the topology used
by Salinas. 
Section 5 deals with certain consequences of the main
result.  We prove that a $KK$-equivalence must
be a homeomorphism. In \S 6 we introduce polonais and
pseudopolonais groups, observe that our results demonstrate that
the $KK$-groups are pseudopolonais, and then state some powerful
results that we learned from C.C. Moore which demonstrate why it
is useful to know the polonais properties.
We use the fact that the $KK$-groups are 
pseudo\-polonais  to demonstrate our primary result,
that the $KK$-pairing is jointly continuous in the natural Salinas
topology.
In \S 7 we conclude by
demonstrating that various structural maps are also continuous,
and finally by showing that the index map is also continuous.

\medbreak\medbreak It is a pleasure to acknowledge  our dependence
upon Salinas's work \cite {19}. Without it this paper would not
exist. We are very grateful    to Larry Brown, Marius Dadarlat,
Gert Pedersen, Chris Phillips, Jonathan Rosenberg, Norberto
Salinas, Bert Schreiber, Klaus Thomsen, and Richard Zekri  for
helpful comments.

In this paper all $C^*$-algebras are assumed separable with the
exception of
those that obviously are not (namely multiplier algebras $\Cal MB$ and
their quotients 
the corona algebras
$\Cal QB = \Cal MB/B$).
If $A$ is not nuclear then when speaking of extensions we always
require that they be semi-split, so that equivalence classes of
extensions coincide with the relevant Kasparov group. All $C^*$-algebras
are assumed to be trivially graded. Tensor products $A\tensor B$ are
always understood to be the minimal tensor product.
 Isomorphisms of topological groups   are
  isomorphisms of groups which are homeomorphisms  as spaces.

\newpage\beginsection {2. Function Spaces }

In this section we review certain aspects of the topologies of the function space
of $*$-homomorphisms for separable $C^*$-algebras $A$ and $B$.

Let $A$ and $B$ be $C^*$-algebras and let $Hom(A,B)$ denote the set of all
 $*$-homo\-morphisms
from $A$ to $B$. There are several possible topologies on $Hom(A,B)$:

\roster
\item The topology of pointwise convergence;
\item The compact-open topology;
\item The topology of uniform convergence on compact sets;
\item If $A$ is separable with countable dense set $\{ a_i \}$
(with each $a_i \neq 0$), the metric topology
    obtained from the metric
$$
\mu (f_1, f_2) = \sum   \frac {| (f_1 - f_2)(a_i)|}{2^i |a_i|}  .
$$
\endroster

\proclaim {Proposition 2.1} The first three topologies on
$Hom(A,B) $ are homeomorphic. If $A$ is separable then all four
topologies are homeomorphic.
\endproclaim

\demo{Proof} Topologies 2), 3), and 4)  agree by standard arguments,
 and of course 3) implies 1), so the only statement requiring proof
is that 1) implies 3): pointwise convergence implies uniform convergence on compact
sets. This is a folklore statement: we insert a proof for convenience.

Suppose that $f_\alpha \to f$ pointwise in $Hom(A,B)$, $K$ is a
compact set, and $\delta > 0$. The set $K$ is a subset of a
separable space, hence separable. Let $\{ a_i \}$ be a countable
dense set for $K$.  Cover $K$ with balls of radius $\delta $
centered at the points $a_i$. Since $K$ is compact, a finite
number of these balls cover $K$, say corresponding to the finite
sequence $a_1, ... , a_n $. Now fix any $a \in K$. Then $|a - a_i
| < \delta $ for some $i$, with $1 \leq i \leq n$. Pick one. Then
$$
\align |(f_\alpha - f)(a) |
&= |(f_\alpha - f)(a - a_i) + (f_\alpha - f)(a_i) |  \\
&\leq  2|a - a_i| +  |  (f_\alpha - f)(a_i)| \\
&\leq 2\delta + |  (f_\alpha - f)(a_i)|   .
\endalign
$$
Now $f_\alpha \to f$  uniformly on the (finite!)
set $\{a_1, \dots , a_n \} $,
and hence $f_\alpha \to f$  uniformly on $K$. \qed\enddemo

Henceforth when we refer to $Hom(A,B)$ as a topological space the
topology described in 1), 2), 3) (or 4) if $A$ is separable) above
is intended and we use their properties interchangably without
further remark. Given a compact set $K \subset A$ and an open set
$U \subset B$ we let

$$ (K,U) = \{ f: A \to B : f(K) \subset U  \} $$

These sets form a subbasis for the topology of $Hom(A,B)$ and in
fact we may require $K$ to be restricted to single points.

Suppose that $f: X \to Y$ is a map of topological spaces, and that
$Im(f)$ denotes the image of the map $f$. There are two obvious
topologies on  this space. One may regard $Im(f)$ as a quotient
space of $X$, in which case the map $X \to Im(f) $ is continuous
and a quotient map, or one may regard $Im(f)$ as a subspace of $Y$
with the relative topology. We refer to these as $Im(f)_{quot}$
and $Im(f)_{rel}$ respectively. The identity map $$ Im(f)_{quot}
\longrightarrow  Im(f)_{rel}
 $$
is of course a continuous bijection. If it is a homeomorphism then
we say that the map $f$ is {\it{relatively open}}.  If $f$ is
actually an open map then it is relatively open, but not
conversely (for instance, a linear  inclusion of a line into the
plane with the standard metric on both is relatively open but not
open.) If $f: X \to Y$ is surjective then open is of course
equivalent to relatively open. If $f: X \to Y$ is injective then
relatively open is equivalent to $f$ being a homeomorphism onto
its image.

Suppose that $f: A \to B$ is a $*$-homomorphism. Then $f(A)$ is
closed in $B$ and hence may be regarded as a $C^*$-algebra in two
ways: either as $A/Ker(f)$ or as $f(A) \subseteq B$. These ways
are algebraically isomorphic and (by the uniqueness of the
topology on a $C^*$-algebra) must be topologically the same as
well. Thus any $*$-homomorphism $f:  A \to B$ must be relatively
open.

\medbreak

\proclaim{Proposition 2.2} Let $A$, $B$, and $C$ be
$C^*$-algebras. Then the composition map $$ Hom(A,B) \times
Hom(B,C) \overset{T}\to\longrightarrow Hom(A,C) $$ is jointly
continuous.
\endproclaim
\medbreak

\demo{Proof} This is essentially identical to \cite {10, p. 259}.
 Let $f : A \to B$, $g : B \to C$, choose some $a\in
A$, and let $W$ be an open neighborhood of $gf(a) \in C$. Then
$g^{-1}W$ is open in $B$ and contains $f(a)$. Choose an open set
$V \in B$ with $f(a) \in V \subset {\overline {V}} \subset g^{-1}W
$. Then $f \in (a,V)$, $g \in (\overline V, W) $ and
$$
 T((x,V),
(\overline V, W)) \subset (a,W)$$ as required.  \enddemo\qed

\proclaim{Proposition   2.3 } Suppose that $h : B \to B'$ is a map
of $C^*$-algebras. Then the induced map $$ h_* : Hom(A,B) \to
Hom(A,B') $$ is continuous. If $h$ is mono then $h_*$ is
relatively open.
\endproclaim

\demo{Proof} Continuity is obvious from the previous proposition,
since joint continuity implies separate continuity. We must show
that $h_*$ is
 relatively open provided that $h$ is mono.  Let
$a \in A$ and let $U$ be an open set in $B$.  Then

$$ \align h_*(a, U ) &=  \{ hf \in Hom(A, B') : f(a) \in U \} \\
&=  \{ hf \in Hom(A, B') : hf(a) \in h(U) \} \quad\text{since $h$
is mono} \\
            &= (a, h(U) ) \cap Im(h_*)  \\
        &= (a,V) \cap Im(h_*)
\endalign
$$
 where $V$ is some open set in $B'$ with
$$
h(U) = Im(h_*) \cap V.
$$
(Such a set $V$ exists since $h$ is a relatively open map, as remarked above.)
Hence $h_*(a,U)  $  is relatively open in $Hom(A,B')$ as required.
\qed\enddemo

\proclaim {Proposition 2.4} Suppose that $h : A' \to A$ is a map
of $C^*$-algebras. Then the induced map $$ h^* : Hom(A,B)
\longrightarrow Hom(A', B) $$ is continuous and relatively open.
\endproclaim

\demo{Proof} Again, continuity is clear from Proposition 2.2. To
show that the map is relatively open,
 let $(a,U) \subset Hom(A,B) $ be a subbasic open set. Then
$$
h^*(a,U) = h^*\{f: A \to B : f(a) \in U  \}
          = \bigcup _x  \,[ (x, U) \cap Im (h^*)  ]
$$
where the union is over all $ x \in h^{-1}(a) $. Each of the sets
$ (x, U) \cap Im (h^*)  $ is relatively open and hence their union
is also relatively open in $Im(h^*)$. \qed\enddemo \medbreak

\proclaim {Proposition 2.5}
 The natural map $$ \Psi
_C : Hom(A,B) \longrightarrow Hom(C\tensor A, C\tensor B ) $$ is
continuous. In particular, the suspension map
$$ S = \Psi _{C_o(\Bbb R)}  : Hom(A,B)
 \longrightarrow Hom(SA, SB)
$$
is continuous.
\endproclaim

\demo{Proof} This is immediate from definitions. \enddemo\qed
\medbreak

\proclaim{Proposition 2.6} Suppose that $A$, $B$, and $C$ are
$C^*$-algebras. Then the composition map $$ [A,B] \times [B,C]
\overset{\overline T}\to\longrightarrow [A,C] $$ is
separately  continuous.
\endproclaim

\medbreak

\demo{Proof}
This is immediate from the definition of the quotient topology.
 \enddemo\qed

Note that we do not claim that $\overline T$ is jointly continuous.
This would be true if 
the map
$$
\CD
 Hom(A,B) \times Hom(B,C)   \\
 @VV{\pi _1 \times \pi _2}V  \\
 [A,B] \times [B,C]
\endCD
$$
 were
a quotient map. In general the product of quotient maps is not necessarily
itself a quotient map. Eventually we will show that the $KK$-pairing
is jointly continuous, but this result will use the fact that
the $KK$-groups are pseudo\-polonais.

\medbreak

 Next we wish to consider the natural map
$$
Hom(A,B) \to Hom(\Cal MA, \Cal MB )
$$
where $\Cal M$ denotes the multiplier algebra of a $C^*$-algebra.
Unfortunately there is no such functor in general, since not every
$C^*$-map from $A$ to $B$ extends to the multiplier algebras.
Fortunately enough maps do, and this leads us to a way to proceed.

The following ideas are due to K. Thomsen \cite{28}, and we are
deeply grateful to him for his assistance. He cites Nigel Higson
\cite{12} as the first to consider such maps.
 A $C^*$-algebra map $h: A \to B$ is said to be
{\it{quasi-unital}} if the closed linear span of $h(A)B$ is
 of the form $pB$ where $p \in \Cal MB$ is some projection.
A quasi-unital map extends to the multiplier algebra level via a
$*$-homomorphism
$$
\Cal Mh : \Cal MA  \to \Cal MB
$$
and hence induces a map at the level of corona algebras denoted
$\Cal Qh: \Cal QC  \to \Cal QD $.

\proclaim {Proposition 2.7} K. Thomsen \cite {28, Prop. 2.8 } If
$A$ and $B$ are $\sigma$-unital and $B$ is stable then any
$C^*$-map $A \to B$
 is homotopic to a quasi-unital map.
If two quasi-unital maps are homotopic then they are homotopic via
a homotopy through quasi-unital maps.
\endproclaim
\qed

Let $Hom(A,B)_{qu} $ denote the quasi-unital maps from $A$ to $B$,
topologized as a subspace of $Hom(A,B)$. Note that $Hom(A,B)_{qu}
$ is {\it{not}} an open subspace of $Hom(A,B)$ in
general\footnote{ There is a counterexample due to K. Thomsen.}
 and hence the inclusion map is not an open map in general.
 We
 let $[A,B] $ denote homotopy classes of $*$-homomorphisms
from $A$ to $B$, topologized as the quotient of $Hom(A,B)$.
Similarly, let
 $[A,B]_{qu} $ denote quasi-unital
homotopy classes of quasi-unital $*$-homomorphisms from $A$ to
$B$, topologized as the quotient of $Hom(A,B)_{qu }$. Then
Proposition 2.7 implies that
 if $B$ is stable then
there is a natural bijection
$$
[A,B]_{qu} \longrightarrow [A,B]
$$
It is obviously continuous. However, it is probably not a
homeomorphism in general,
 by the above remarks.

\medbreak

\proclaim{Proposition 2.8} Suppose that $f: A \to B $ is
quasi-unital. Then so is the map $1\tensor f: SA \to SB$.
Conversely, if $f: A \to B$ is a $*$-homomorphism and $1\tensor f$
is quasi-unital then so is $f$.
\endproclaim

\demo{Proof} I am endebted to Klaus Thomsen for the following
proof.

If $f$ is quasi-unital there is some projection $p \in \Cal MB$
with $$ \overline {f(A)B }   =   pB $$ where $\overline {f(A)B } $
denotes the closed linear span of $f(A)B$. Then a direct
calculation shows that $$ \overline{ (1\tensor f)(SA)SB } =
(1\tensor p)SB $$ and $1\tensor p \in   \Cal M(SB)$.

In the other direction, suppose that $1\tensor f$ is quasi-unital
with
$$
\overline{f(SA)SB} = eSB
$$
 for some projection $e \in \Cal M (SB)$.
Let
$$
q : SB = C_o(\Bbb R)\tensor B \to B
$$
be evaluation at zero. Then $q$ is surjective, and
$$
{f(A)B} = q(1\tensor f)(SA)(SB)
$$
 as required. \qed\enddemo \medbreak

\medbreak

{\bf{Remark 2.9}}  Quasi-unital maps give a nice insight into the
general theory of extensions of $C^*$-algebras. An essential
extension
$$
\tau :   \qquad\qquad  0  \to B\tensor\Cal K  \to E  \to  A  \to
0
$$
with Busby classifying map
$$
\tau : A \longrightarrow \Cal Q(B\tensor\Cal K)
$$
 lies in the group $\KKgraded 1.A.B $.
If $f:B \to B'$ is a $*$-homomorphism then $f_*[ \tau ] \in \KKgraded
1.A.{B'} $. However, it is not so easy to see directly how to
construct the corresponding extension or its classifying map
$$
\hat\tau : A \longrightarrow \Cal Q(B'\tensor\Cal K)  .
$$
The natural thing to do is to take a pushout construction, as one
would do in abelian groups, but the diagram
$$
\CD
 0  @>>>  B\tensor\Cal K   @>>>   E   @>>>   A   @>>>   0  \\
 @. @VV{f\tensor 1}V    \\
@. B'\tensor\Cal K
 \endCD
 $$
 may not be completed without some assumptions on $f$. Eilers, Loring,
  and Pedersen \cite {11}
  show that if $f$ is {\it{proper}}
 then it is possible to complete the pushout
 and thereby obtain $\hat\tau $.
    However, not every
 map is proper.

 Using our knowledge of quasi-unital maps, there is a nice
 solution.  Given $f$ as above, replace it by a quasi-unital
 map $g$ homotopic to it. This then extends to
 $$
 \Cal Mg : \Cal
 M(B\tensor\Cal K) \to \Cal M(B'\tensor\Cal K)
 $$
  and hence
 induces a map
$$
 \Cal Qg : \Cal
 Q(B\tensor\Cal K) \to \Cal Q(B'\tensor\Cal K)
 $$
Then we have $$ A \overset\tau\to\longrightarrow Q(B\tensor\Cal K)
\overset{\Cal Qg}\to\longrightarrow \Cal Q(B'\tensor\Cal K)
$$
If this map happens to be mono then it classifies an extension
which is in the class $f_*[\tau ]$. If not (as is likely), simply
add on a map of the form
$$
A \overset{\pi\sigma}\to\longrightarrow \Cal Q(B'\tensor\Cal K)
$$
where
$$
\sigma : A \to \Cal M(B'\tensor\Cal K )
$$
is some trivial extension.  Then the composite
$$
(\Cal Qg)\tau \oplus \pi\sigma : A \longrightarrow \Cal
Q(B'\tensor\Cal K) \oplus    \Cal Q(B'\tensor\Cal K) \to \Cal
Q(B'\tensor\Cal K)
$$
is mono and it represents the desired element:
$$
f_*[\tau ] = [(\Cal Qg)\tau \oplus \pi\sigma ] \in \KKgraded 1.A.{B'}   .
$$

\newpage\beginsection{ 3. Continuity of the Cuntz pairing}

In this section we introduce the Cuntz quasi-homomorphism picture
of the Kasparov groups and verify that the $KK$-pairing is
continuous in the corresponding topology.

Let $A$ be a $C^*$-algebra and let $A\star A$ denote the
$C^*$-algebra free product. Then there are two canonical maps
$\alpha _A $ and $\bar\alpha _A $ mapping $A \to A\star A$. (When
the context is clear we omit the subscript.)  Following J.Cuntz
\cite{6, 7, 8},  let $qA $ denote the closed ideal of $A\star
A$ generated by all elements of the form $\alpha (a) -  \bar\alpha
(a) $. The identity map $A \to A$ induces a canonical map $A\star
A \to A$ and a natural short exact sequence $$
 0  \to qA \to A\star A \to A \to 0.
  $$
   It
is easy to see that $q$ is a covariant functor and that if $A$ is
separable then so too is $qA$.

Cuntz has shown that there is a natural isomorphism

 $$
   [qA, B\tensor\Cal K]   \overset\cong\to\longrightarrow
     \KKgraded 0.A.B     .
\tag 3.1
 $$
We give $\KKgraded 0.A.B $ a topology by declaring the map (3.1)
to be a homeomorphism. We shall refer to this topology as the
{\it{Cuntz topology }} on $\KKgraded 0.A.B $. Similarly we
topologize $\KKgraded 1.A.B $ by declaring the isomorphism
$$
\KKgraded 1.A.B   \cong  \KKgraded 0.SA.B  
 \cong [ q(SA), B\tensor\Cal K ]
$$
to be a homeomorphism. (In Section 4 we will discuss other
possibilities for the topology on $\KKgraded *.A.B $ and at that
time we will refer to the present topology as $\KKgraded *.A.B _q
$.)

Next we need some elementary facts.
 \proclaim{Proposition 3.2} Suppose that $A$ is separable. Then
   the
natural map
$$
 q: Hom (A,B) \to Hom
(qA, qB)
$$
 and the induced map
  $$
 q_*:   [A,B] \to [qA, qB]
   $$
    are both continuous.
\endproclaim

\medbreak

\demo{Proof}  It is clear that the natural map
$$
Hom(A,B) \longrightarrow Hom(A\star A, B\star B )
$$
given by $f \to f\star f $ is continuous.
Composing with the natural map $qA \to A\star A $
 yields
a continuous map
$$
Hom(A,B) \longrightarrow Hom(qA, B\star B)
$$
which in fact factors through $Hom(qA, qB)$. As
$Hom(qA, qB) $ has the relative topology in $Hom(qA, B\star B)$,
this implies that the natural map
$q : Hom(A,B) \to Hom (qA, qB) $ is continuous.

Since the quotient maps $  Hom(A,B) \to [A, B] $ and 
$ Hom(qA,qB)
\to [qA, qB] $ are continuous,   it
follows immediately
 that the induced map  $q_*: [A,B] \to [qA, qB] $ is
also  continuous.
 \enddemo \qed

\medbreak

Next we introduce a key map

$$
 \Phi _{AB} : [qA, B\tensor\Cal K ] \longrightarrow [qA, qB\tensor\Cal K] .
 $$
 \medbreak\flushpar
  We
follow the notation of \cite{13, Theorem 5.1.12}.  There are
canonical (at least up to unitary equivalence) $*$-homomorphisms
$$ \rho _B : q(B\tensor\Cal K) \to qB\tensor\Cal K, \qquad \phi _B:
qB \to q^2B \tensor  M_2(\Bbb C) , \qquad \theta : \Cal K \tensor
M_2(\Bbb C) \to \Cal K . $$
 Using these, we
define $\Phi _{AB}$ to be the composite

$$ \align [qA, B\tensor\Cal K ] &\overset{q}\to\longrightarrow
[q^2A , q(B\tensor\Cal K)]  \\
 &\overset{(\rho _B)_*}\to\longrightarrow [q^2A, qB \tensor
\Cal K ]   \\ \
 &\overset{\Psi_{M_2(\Bbb C)}}\to\longrightarrow [q^2A\tensor
M_2(\Bbb C), B\tensor\Cal K \tensor M_2(\Bbb C) ]   \\
 &\overset{(\phi
_A)^*}\to\longrightarrow [qA, B\tensor\Cal K\tensor M_2(\Bbb C)
]\\
 &\overset{\theta
_*}\to\longrightarrow  [qA, B\tensor\Cal K ] .
\endalign
 $$

It is immediate from Proposition 3.2 and results of \S 2 that
$\Phi _{AB} $ is continuous. In fact more is true:

\proclaim{Proposition 3.3} The natural map  $$
 \Phi _{AB} : [qA, B\tensor\Cal K ] \longrightarrow [qA, qB\tensor\Cal K] .
 $$
 is an isomorphism of topological groups.
 \endproclaim

 \demo{Proof} We have shown in the discussion above that the map
$\Phi _{AB} $  is continuous.
 There is a natural inverse map
which is constructed as follows.
  The natural maps $1_B : B \to B $
 and $0 : B \to B$ induce a natural map $1_B \star 0 : B\star B \to B$.
 Let
 $$
 \gamma ^B : qB \to B
 $$
 be the restriction of this map to $qB$.  Then $\gamma ^B$
 induces a continous map
 $$
 (\gamma ^B \tensor 1)_* :
  [qA , qB\tensor\Cal K ] \longrightarrow [qA, B\tensor\Cal K
 ]
$$ which is shown in \cite {13}   to be the map inverse to $\Phi _{AB} $.
The map $(\gamma ^B \tensor 1)_*$ 
is continuous by Propositions 2.3 and  2.5.
\enddemo\qed

\medbreak

\proclaim{Theorem 3.4} The Kasparov pairing $$ \KKgraded *.A.B
\,\,\times\,\, \KKgraded *.B.C \overset{\tensor
_B}\to\longrightarrow \KKgraded *.A.C
$$ is separately continuous in the Cuntz topology.
\endproclaim
\medbreak

\demo{Proof} Suppose first that $* = 0$ so that we are looking at
the
 pairing
 $$
  [qA, B\tensor\Cal K ] \,\,\times\,\,[qB,
C\tensor\Cal K ] \longrightarrow [qA, C\tensor\Cal K ].
 $$
This pairing is simply the composite \medbreak

 $$
  \CD [qA , B\tensor\Cal K
]\,\,\times\,\, [qB, C\tensor\Cal K ]
\\
@VV{\Phi _{AB} \tensor \Psi _{\Cal K}}V   \\
  [qA ,qB\tensor\Cal K ]\,\,\times\,\, [qB\tensor\Cal K, C\tensor\Cal
K\tensor\Cal K  ]  \\
 @VV{compose}V   \\
[qA,C\tensor\Cal K\tensor\Cal K ]    \\
 @VV{\cong}V  \\
 [qA,
C\tensor\Cal K ]
\endCD $$
\medbreak\medbreak\flushpar
 and hence is continuous in each variable. Since Bott periodicity (in
either variable) may be regarded as pairing with $KK$-elements of
degree zero, it follows at once that the natural Bott map
 $$
  \KKgraded 0.A.B  \to  \KKgraded 0.A.{S^2B}
$$
 is a homeomorphism, and
similarly in the $A$ variable.

Next, we recall that the Cuntz topology on $\KKgraded 1.A.B$ was
given
 by
insisting that the natural connecting  isomorphism
$$
\KKgraded 0.{SA}.{B}       \overset\cong\to\longrightarrow
 \KKgraded 1.A.B
$$
 which arises from the canonical short exact
sequence
 $$ 0 \to SA \to CA \to A \to 0
 $$
 be a homeomorphism. (This choice is of course consistent with the
 Bott maps.) Then it is an easy exercise to prove that the pairing
 is continuous as stated.
 \enddemo\qed

The most general case of the Kasparov pairing is built from the
previous pairing and the following map.

\proclaim{Proposition 3.5} The external mapping $$
\tau _C : \KKgraded *.A.B
\to \KKgraded *.{A\tensor C}.{B\tensor C}  $$ is continuous.

\endproclaim

\medbreak \demo{Proof} There is a natural map $\nu : q(A\tensor C)
\to qA \tensor C $, which arises as the restriction of the natural
map $$ (\alpha _A \tensor 1_C) \star (\bar\alpha _A \tensor 1_C) :
(A\tensor C)\star (A\tensor C)  \to  (A\star A)\tensor C $$ to
$qA$
 and the mapping above is the composite
 $$
[qA,B\tensor\Cal K ] \overset{\Psi_C}\to\longrightarrow [qA\tensor
C, B\tensor C\tensor\Cal K ]  \overset{\nu^*}\to\longrightarrow
[q(A\tensor C), B\tensor C\tensor\Cal K ] $$ hence continuous.
\enddemo\qed

\medbreak

Combining results we have the most general theorem of this section.

\medbreak

\proclaim {Theorem 3.6} The Kasparov pairing
$$
\KKgraded *.{A_1}.{B_1\tensor D} \times \KKgraded *.{D\tensor A_2}.{B_2}
\overset{\tensor _D}\to\longrightarrow
\KKgraded *.{A_1\tensor A_2}.{B_1\tensor B_2}
$$
is separately continuous in the topology associated to the Cuntz
picture of
$$
\KKgraded 0.A.B \cong [qA, B\tensor\Cal K ].
$$
\endproclaim

\medbreak
\demo {Proof} Let $E = B_1\tensor D \tensor A_2$. Then
$$
x \tensor _D y = \tau _{A_2}(x) \tensor _E \tau _{B_1}(y)
$$
so that 3.4 and 3.5 together imply the result.
\enddemo\qed

\medbreak
We shall improve this theorem to show joint continuity in Theorem 6.8.

\newpage\beginsection{4. The
Topologies on $\KKgraded 1.A.B $}

\medbreak

In this section we show that four topologies on $\KKgraded 1.A.B $
coincide. Of course the topology that we are interested in is the
topology used by Salinas, which we shall denote $\KKgraded 1.A.B
_S $ in this section. The main consequence of this section is that
this topology coincides with the natural topology on $[q(SA),
B\tensor\Cal K ]$ in the Cuntz picture of $KK$. Since we have
already shown that the Kasparov pairing is continuous in the Cuntz
picture, this will imply that the pairing is continuous in the
Salinas topology. This is what we need in applications.

Here are the four topologies that we shall consider:

\medbreak

{\bf 1. The Salinas Topology:\quad $\KKgraded 1.A.B  _S $ }

\medbreak

We regard $\KKgraded 1.A.B $ as the quotient of the space
 $\Cal E(A,B) \subseteq Hom(A,\Cal Q(B\tensor \Cal K))$
of associated semisplit extensions, with natural
 quotient map $$ \Upsilon _S: \Cal
E(A,B)  \longrightarrow \KKgraded 1.A.B       . $$ We  topologize
$\Cal E(A,B) $ by giving it the relative topology in
$Hom(A,\Cal Q(B\tensor\Cal K))$ or equivalently by giving it
a metric topology as in (2.1).
Fix one such metric $\mu $ on $\Cal
E(A,B)$.
 Given $x$, $x' \in \KKgraded 1.A.B   $, define
$$ \hat \mu (x,x') = \,{\text{inf}} \,\mu (\tau , \tau ') $$ where
the $inf$ is taken over all $\tau \in x$, $\tau ' \in x'$. Salinas
shows \cite {19, 3.1} that one obtains a pseudo\-metric\footnote{
A {\it {pseudo\-metric}} is a  function $\mu : X\times X \to \Bbb
R$ satisfying the usual metric axioms except that if $\mu (x,y) =
0$ it need not be the case that $x = y$. If $X$ is a topological
group we insist that the pseudo\-metric be continuous and that it
be invariant under the group action: $\mu (xz, yz) = \mu (x,y)$.}
 on
$\KKgraded 1.A.B $. The associated  topology  is of course
 independent of
the choice of
 the sequence  $\{a_i\}$.   With
  respect to the pseudo\-metric   $\KKgraded 1.A.B $
is a topological group, and the pseudo\-metric is invariant under
the group action. The quotient topology induced upon $\KKgraded
1.A.B $ from $\Cal E(A,B) $
 via $\Upsilon _S $
coincides with this topology. We shall denote $\KKgraded 1.A.B $
with this topology as $\KKgraded 1.A.B  _S$.

Salinas shows that this construction
defines a functor to the category of topological groups and
continuous homomorphisms in each variable. One special case is
easy to describe: the group $K_1(B) = \KKgraded 1.{\Bbb C}.B _S$ is
countable, complete, and hence discrete. However in general the
topology may be highly non-trivial.

We note that in Salinas's paper, he has a standing assumption that
$A$ must be unital. Dadarlat (\cite{9} and private communication)
has verified that this assumption is not essential. All of
Salinas's results remain true in the nonunital case, including the
equality between the closure of zero and the quasidiagonal
extensions (cf. \cite {25}). \medbreak

{\bf 2. The Zekri topology \quad $\KKgraded 1.A.B  _Z $}

\medbreak

Following Zekri \cite{29, 30},  let $EA$ denote the universal
$C^*$-algebra generated by a separable $C^*$-algebra $A$ and an
element $F$ with $F^* =F$ and $F^2 = 1$, with $1$ acting on $A$ as
the identity. Let $[A,F]$ denote the set of commutators and let
$\epsilon A$ denote the smallest closed ideal of $EA$ which
contains $[A,F]$. Then there is a canonical commutative diagram
 $$
  \CD 0  @>>>   \epsilon A
@>>>  EA  @>>> A^+\oplus A^+  @>>>  0
\\ @.  @AAA    @AAA @AAdA   \\
  0  @>>>   \epsilon A  @>>>  \tilde EA  @>>> A  @>>>
0   \\
\endCD
$$
 where $d(a) = (a,a) $ is the diagonal map,
 $A^+$ is the unitalization of $A$,  and $\tilde EA$ is the
pullback.

Zekri  defines a {\it{K-extension }} $e$  of $B$ by $A$
 to be a diagram
$$
 \CD 0 @>>>   A   @>>>   C^*(A,U)   @>>>
A^+\oplus A^+   @>>>   0   \\ @.  @.   @AAA   \\ @.  @. J   @>{\mu
}>>   B\tensor\Cal K
\endCD
\tag 4.1
$$
where $U = U^* = U^{-1}$ is a self-adjoint unitary, $J$ is the
ideal generated by $[A,U]$, the row is exact, and $\mu $ is a
monomorphism. Diagram (4.1) yields a universal example of the form
$$
 \CD 0 @>>>   A   @>>>   C^*(A,U)   @>>>
A^+\oplus A^+   @>>>   0   \\ @.  @.   @AAA   \\ @.  @.  \epsilon
A   @>{1}>>   \epsilon A
\endCD
$$
(after stabilizing $A$). Let $KExt(A,B)$ denote
the set of  K-exten\-sions.

There is a natural map
 $$
 \zeta : KExt(A,B) \longrightarrow
Hom(\epsilon A, B\tensor \Cal K )
 $$
  defined as follows. Given a
K-extension $e$
as in (4.1),
there is a canonical map $EA \to C^*(A,U) $. Its
restriction to $\epsilon A$ factors uniquely through $J$. Then we
define $\zeta (e)$ to be the composite
$$
 \epsilon A  \to J \to
B\tensor \Cal K .
 $$
\medbreak

\medbreak \proclaim{Proposition 4.2 } (Zekri \cite {29, 
Theorem 2.4}) The natural map
$$
 \zeta : KExt(A,B) \longrightarrow
Hom(\epsilon A, B\tensor \Cal K )
 $$
is a bijection on isomorphism classes. \qed\endproclaim

Zekri shows 
\cite {29, Theorem 4.4} that if $A$ is $K$-nuclear 
then there is a natural isomorphism
$$
[\epsilon A, B\tensor\Cal K ] \cong \KKgraded 1.A.B
$$
and since the left hand side has a natural topology as a quotient
of
the space
$Hom (\epsilon A, B\tensor\Cal K ]$, the right hand side
inherits a topology via this isomorphism which we denote
$\KKgraded 1.A.B _Z $.

\medbreak

{\bf 3. Another Cuntz topology:  \quad $\KKgraded 1.A.B  _C $}
\medbreak

This topology arises from the  Cuntz picture of $\KKgraded 1.A.B
$. (It is a priori different from the topology that we discussed
in \S 3.)

Let $H_B$ denote the universal Hilbert $B$-module, and let
$$
P = \left[\smallmatrix 1&0 \\ 0 & 0 \endsmallmatrix\right] \in
\Cal M(H_B \oplus H_B) .
$$

Let $\Bbb
E_1(A,B)$ denote the collection of maps
$$
 \sigma : A \to \Cal M(H_B \oplus H_B )
 $$
such that the commutator
 $[\sigma (A), P] \subset B\tensor \Cal K$, and
 the resulting $*$-homomorphism
 $$
 \Upsilon _C (\sigma ) : A \to \Cal Q(B\tensor K)
 $$
given by
$$
\Upsilon _C (\sigma)(a) = \pi(P\sigma (a) P)
$$
 is a monomorphism.
 This gives a natural
map
$$
 \Upsilon _C :      \Bbb E_1(A,B)  \to  \KKgraded 1.A.B
$$
 We topologize the set $\Bbb E_1(A,B)$ as follows. Let 
 $\bar\mu $ be a standard metric on the space $Hom(A, \Cal M(H_B))$
 as in (2.1). Define a pseudometric on $\Bbb E_1(A,B)$ by
 $$
 \mu ((\sigma , P),(\sigma ' , P)) = 
 max \{ \bar\mu (P\sigma P , P\sigma ' P),\,\,
 \bar\mu (   (1 - P)\sigma (1 - P) , (1 - P)\sigma ' (1 -P) )
 \} .
 $$
 \medbreak\flushpar
Then $\KKgraded 1.A.B $ inherits a topology as the quotient of
$\Bbb E_1(A,B) $ which we denote $\KKgraded 1.A.B  _C $.

\medbreak

{\bf 4. The $qA$ topology: \quad $\KKgraded 1.A.B  _q $}

\medbreak

This is the topology that we used in Section 3, for which the
$KK$-product is separately continuous.
Recall that
this topology is simply described. Cuntz has  shown that
there is a natural isomorphism $[qA, B\tensor\Cal K] \cong \KKgraded 0.A.B $
and hence a natural isomorphism
$$
\KKgraded 1.A.B  \cong \KKgraded 0.{SA}.B  \cong [q(SA), B\tensor\Cal K ] .
$$
We denote the resulting topology by $\KKgraded 1.A.B _q $.

\medbreak

We must show that these four topologies coincide.
Here is the first step:

\medbreak

\proclaim{Proposition 4.3} There is a  canonical map
$$
h_* :  \Bbb E_1(A,B) \longrightarrow \Cal E(A,B)
$$
which is continuous and which induces an isomorphism
of topological groups
$$
\tilde h :  \KKgraded 1.A.B  _C  \longrightarrow \KKgraded 1.A.B
_{S}.
$$
\endproclaim

\medbreak

If we ignore topologies, then $\tilde h$ is just the identity map.
\medbreak

\demo{Proof}  There is a natural map $\Cal M(H_B \oplus H_B ) \to
\Cal M(H_B) $ given by
$$
T = \bmatrix T_{11}   &   T_{12}  \\  T_{21}  &  T_{22}
\endbmatrix   \longrightarrow   PTP  =  T_{11} .
$$
This is not a $*$-homomorphism, of course, but it certainly is
continuous. Composing with the projection $\Cal M(H_B) \to
\Cal Q(H_B)$
we obtain a map
$$
h: \Cal M(H_B \oplus H_B )  \to \Cal Q(H_B)
$$
 which is continuous. Then $h_*$ is defined\footnote{Although
  we do not need this, we note that 
 $h_*$ is onto.    Any semisplit $C^*$-monomorphism
$$
\tau : A \to \Cal Q(H_B) 
$$
 lifts to a completely positive contraction $A
\to \Cal M(H_B) $. Then this completely positive map may be
dilated to a $C^*$-homomorphism 
$$
\sigma : A \to \Cal M(H_B\oplus H_B) 
$$
 by
the generalized Stinespring theorem, due to Kasparov \cite{14}.
Then $(\sigma , P) \in \Bbb E(A,B)$ and $h_*(\sigma , P) = \tau $.  
 Thus the canonical map
$$
h_* : \Bbb E_1(A,B) \longrightarrow \Cal E(A,B)
$$
is continuous and surjective,} by
 $$
 h_*(\sigma , P) = h\sigma  
 $$
and it induces a continuous map
$$
\tilde h :  \KKgraded 1.A.B  _C  \longrightarrow \KKgraded 1.A.B
_S .
$$
Of course this is the identity map at the level of sets, but we
must keep track of the topologies as well.
 To prove that $  \tilde h $
 is an isomorphism of topological groups, it suffices then to
prove that the map $\tilde h$ is open. By homogeneity, it suffices
to show this in a neighborhood of zero.

So fix $\epsilon > 0$ and let $B_\epsilon ^C $ be an open ball
about zero in $\KKgraded 1.A.B  _C $. We claim that
$$
\tilde h (B_\epsilon ^C  ) =         B_\epsilon ^S
$$
where $ B_\epsilon ^S$ is the corresponding open ball in
$\KKgraded 1.A.B _S$. (This of course implies that $\tilde h$ is
open.)

Suppose first that $[\tau ] \in B_\epsilon ^C $. Then there is
some trivial extension
$$
(\sigma ,P) \in \Bbb E(A,B)
$$
 and some extension
$$
(\hat \tau , P) \in \Bbb E(A,B)
$$
 with $\pi (P\hat\tau P) = \tau
$ and 
$$
\mu ((\sigma ,P), (\hat\tau ,P) ) < \epsilon   .
$$
Then $\bar\mu (P\sigma P, P\hat\tau P ) < \epsilon $ in 
$Hom(A, \Cal M(H_B))$ and hence $\mu ([P\sigma P], [\tau ]) < \epsilon $.
Thus  $\tilde h[\tau ]\in B_\epsilon ^S  $. So 
$\tilde h(B_\epsilon ^C) \subseteq B_\epsilon ^S $.

Conversely, suppose that  $[\tau ]\in B_\epsilon ^S  $. Then also
$-[\tau ]\in B_\epsilon ^S  $, since the pseudo\-metric is invariant
the group action.
By abuse of language, we write
$-\tau $ for some representative of $-[ \tau ]$. There is some
trivial extension $\sigma : A \to \Cal M(H_B) $ with $\mu (\tau,
\pi\sigma ) < \epsilon $, and similarly, some trivial extension
$\sigma ': A \to \Cal M(H_B) $ with $\mu (-\tau, \pi\sigma ' ) <
\epsilon $. Let
 $$ \bar\sigma = \bmatrix  \sigma   &   0 \\   0  &
\sigma '
\endbmatrix   .
$$

The extension $\tau \oplus (-\tau )$ is a trivial extension and thus
lifts to some
$$
\tilde \tau : A \to \Cal M(H_B \oplus H_B).
$$
 Then $(\tilde \tau ,P) \in \Bbb E(A,B) $. 
 Finally, it is easy to see that
$$
\mu ( (\tilde\tau , P), (\bar\sigma , P)) < \epsilon 
$$
with respect to the pseudometric on $\Bbb E(A,B)$
so that
 $$
\tilde h[\tilde\tau ] =     [\tau ] \in B_\epsilon ^C
 $$
as required.
This implies that $\tilde h(B_\epsilon ^C) = B_\epsilon ^S $
and hence that
$\tilde h$ is  open. As $\tilde h$ is known
to be a continuous
bijection, it must be an isomorphism of topological
groups, as desired. \qed\enddemo

\medbreak

Next we compare the Cuntz and Zekri topologies, and for this we
assume that $A$ is $K$-nuclear.

Define
$$
 \hat\alpha : \Bbb E_1(A,B) \to KExt(A,B)
 $$
  as
follows.  Given $\sigma  \in \Bbb E_1(A,B)$, there is an
associated K-extension
$$
 \CD
 0 @>>> A  @>>> C^*(A, 2P - 1)  \\
 @.  @.   @AAA   \\
   @. @.   B\tensor\Cal K   @>1>> B\tensor \Cal K
\endCD
$$
 which we denote $\hat\alpha (\sigma , P)$.

  The composite
  $$
   \Bbb E_1(A,B)
\overset{\hat\alpha}\to\longrightarrow KExt(A,B)
\overset{\tilde\pi}\to\longrightarrow  [\epsilon A,
B\tensor\Cal K ]
$$
 is given as follows. Given $(\sigma , P) $,
the universal map
$$
 EA \to C^*(A, 2P - 1)
  $$
   induces a
commutative diagram
$$
\CD EA  @>>>  C^*(A, 2P - 1)  \\
@AAA  @AAA   \\
\epsilon A  @>\eta >>   B\tensor K @>1>>  B\tensor K
\endCD
$$
and the homotopy class of $\eta $  is defined to be
 $\alpha (\sigma , P)  \in [\epsilon A , B\tensor\Cal K ]$.

The map $\tilde\pi\tilde\alpha $ is continuous, for if $\sigma _j
$ converges to   $\sigma $ then it is clear that the
associated maps
$$
 \hat\alpha (e_j) : \epsilon A  \to   B\tensor K
  $$
 converge to $\hat\alpha (e)$. Further, the map descends to
 homotopy classes, giving a commuting diagram
$$
 \CD
  \Bbb E_1(A,B)
@>{\hat{\alpha }}>>  KExt(A,B)
 \\ @VV{\pi }V    @VV{\tilde\pi }V  \\
 \KKgraded 1.A.B _C  @>\alpha >> [\epsilon A, B\tensor \Cal K] @>\cong>>
 \KKgraded 1.A.B _Z
\endCD
$$

 We may describe the map $\alpha $ in a more universal fashion.
 Recall that there is a universal K-extension which we denote $\hat
 e$. It is immediate that $\hat\alpha (\hat e) : \epsilon A
 \to \epsilon A $ is the identity element $1_{\epsilon A}$. Then we have
 $$
 \alpha
 (e) = \hat\alpha (e)_* 1 _{\epsilon A}
 $$
  for any extension $A$, where
  $$
\hat\alpha (e)_*  : [\epsilon A , \epsilon A \tensor\Cal K ] \to
[\epsilon A, \epsilon B\tensor\Cal K ].
$$
It is clear from this description that $\alpha $ is continuous.

 Next we define
 $$
\beta : \KKgraded 1.A.B _Z
 \to  \KKgraded 1.A.B _C     .
  $$
  
Recall that the universal $K$-extension has the form 
   $$
   \CD 0  @>>>   A   @>>>   EA
 \\ @.   @.   @AAA  \\ @.   @.   \epsilon A   @>1>> \epsilon A.
 \endCD
 $$
  Expand the diagram to 
  $$
   \CD 0  @>>>   A   @>>>   EA   @>>> \Cal M(\epsilon A)
 \\ @.   @.   @AAA  \\ @.   @.   \epsilon A   @>1>> \epsilon A.
 \endCD
 $$
 Let $\sigma _u : A \to \Cal M(\epsilon A) $ be the composite
 $$
 A \to EA  \to \Cal M(\epsilon A)
 $$
 in the diagram, and let $P_u = \sigma _u ((F + 1)/2 )$. 
 Then $(\sigma _u, P_u)$ defines a universal 
  element  $u_A \in \Bbb E_1(A, \epsilon A)$ and a corresponding
  class $[u_A] \in \KKgraded 1.A.{\epsilon A} _C$. 
  If $h : \epsilon A \to B\tensor\Cal K$ then we define
  $$
  \beta (h) = h_*([u_A]).
  $$

  Alternately we may construct this class explicitly using 
  quasi-unital approximations.
   Given $[h] \in [\epsilon A, B\tensor\Cal K]$, choose a 
   quasi-unital representative $h : \epsilon A \to B\tensor \Cal K$. 
   Then we may define $\beta (h)  = [(\sigma _h , P_h)]$ where
   $\sigma _h$ is the composite
   $$
   A \overset{\sigma _u}\to\longrightarrow \Cal M(\epsilon A)
   \overset{\Cal M(h)}\to\longrightarrow \Cal M(B\tensor \Cal K)
   $$
   and $P_h = \Cal M(h)(P_u) .
   $
   
It is clear from the definitions
 that $\beta $ is continuous.

\medbreak

 \proclaim {Theorem 4.4 } Suppose that $A$ is $K$-nuclear. 
 Then the map
  $$ \alpha : \KKgraded 1.A.B  _Z  \cong [\epsilon A, B] \to
 \KKgraded 1.A.B _C
 $$
  is an isomorphism of topological groups with
 inverse $\beta $.
 \endproclaim
 \medbreak

 \demo{Proof} Zekri has shown that the map $\alpha $ is a bijection
 with inverse $\beta $, so the only point at issue is whether these
 maps are continuous. They are indeed, as we have shown above.
 \enddemo\qed

\medbreak We note that Higson has shown (cf. \cite{2, 22.1}) that
for any choice of a rank 1 projection $e$, the natural map $A \to
A\tensor\Cal K$ sending $a$ to $a\tensor e$ gives rise to a
natural isomorphism
$$
[A \tensor \Cal K , B\tensor\Cal K ]
\overset\cong\to\longrightarrow [A , B\tensor\Cal K ] .
$$
We use this isomorphism below without further comment.
\medbreak

Next we show that $\KKgraded 1.A.B _q \cong \KKgraded 1.A.B _Z $.

\medbreak \medbreak

\proclaim{Theorem 4.5} Suppose that $A$ is $K$-nuclear. Then
 there is a natural homeomorphism
$$
[\epsilon A, B\tensor\Cal K]  \overset\cong\to\longrightarrow [qSA,
B\tensor\Cal K ]
$$
such that the diagram \medbreak
$$
\CD
 [\epsilon A, B\tensor\Cal K]   @>>>  [qSA, B\tensor\Cal K ]  \\
 @VV\cong V        @VV\cong V    \\
 \KKgraded 1.A.B _Z  @>>>    \KKgraded 1.A.B _q
 \endCD
 $$
 commutes.
  \endproclaim

The $K$-nuclearity assumption
(or something like it)
apparently is needed since under
that hypothesis there is a natural homotopy equivalence $\epsilon
A \to q\epsilon A$ which we use.  This is due to Zekri \cite {29},
who notes that the $\epsilon A$ is not necessarily homotopy
equivalent to $q\epsilon A$  for $A$ not $K$-nuclear.

  \medbreak
  \demo{Proof}
Zekri shows \cite{30} that $\epsilon A $ is $KK$-equivalent to
$SA$. Let
$$
u : q\epsilon A  \to  SA \tensor\Cal K    \qquad\text{and}\qquad
v: qSA \to \epsilon A \tensor\Cal K
$$
represent classes $[u] \in \KKgraded 0.{\epsilon A}.{SA} $ and
$[v] \in \KKgraded 0.{SA}.{\epsilon A} $  in the Cuntz picture
which implement the $KK$-equivalence.

Define a map
$$
\eta : [\epsilon A, B\tensor\Cal K] \longrightarrow
[qSA, B\tensor\Cal K]
$$
to be the composite
$$
[\epsilon A  , B\tensor\Cal K] \,\cong\, [ \epsilon A \tensor\Cal
K, B\tensor\Cal K] \overset{v^*}\to\longrightarrow
 [qSA,
B\tensor\Cal K]
$$

In the other direction, define a map
$$
\theta :   [qSA, B\tensor\Cal K]
\longrightarrow  [\epsilon A, B\tensor\Cal K]
$$
as follows. If $g : qSA \to B\tensor\Cal K$
then define $\theta ([g]) $ to be
the homotopy class of the composition
$$
\epsilon A \overset\simeq\to\longrightarrow q\epsilon A
\overset{qu}\to\longrightarrow q(SA \tensor \Cal K)
\longrightarrow q(SA) \tensor \Cal K
\overset{g\tensor 1}\to\longrightarrow B\tensor\Cal K\tensor\Cal K
\to B\tensor\Cal K
$$
where $\epsilon A \overset\simeq\to\longrightarrow q\epsilon A$
is the homotopy equivalence of Zekri. (This is where we use the
$K$-nuclearity assumption.)
These maps are clearly homomorphisms which are natural in $B$ for
fixed $A$. (They do depend upon the choice of $KK$-equivalence, of
course.) We must show that $\eta $ and $\theta $ are inverse to
one another.

The map $\eta\theta ([g])$ is the composite
$$
qSA \overset{v}\to\longrightarrow \epsilon A \tensor\Cal K\simeq
q(\epsilon A) \tensor\Cal K \simeq q^2(\epsilon A) \tensor \Cal K
\overset{qu \tensor 1 }\to\longrightarrow qSA \tensor\Cal
K\overset{g}\to\longrightarrow B\tensor\Cal K \tensor\Cal K \cong
B\tensor\Cal K
$$
which is homotopic to $ qSA \overset{g}\to\longrightarrow
B\tensor\Cal K $ since the composite

$$
qSA \overset{v}\to\longrightarrow \epsilon A \tensor\Cal K\simeq
q(\epsilon A) \tensor\Cal K \simeq q^2(\epsilon A) \tensor \Cal K
\overset{qu \tensor 1 }\to\longrightarrow qSA \tensor\Cal
K
$$
is essentially the $KK$-product
$$
v\tensor _{\epsilon A} u = 1_{\KKgraded 0.{SA}.{SA}  } .
$$
Thus $\Psi\Theta ([g]) = [g] $.
Similarly, in the other direction  $\theta\eta ([f]) = [f]$
via the fact that
$$
u \tensor _{SA} v = 1_{\KKgraded 0.{\epsilon A}.{\epsilon A} }  .
$$
The maps $\eta $ and $\theta $ are visibly continuous, and hence
they are homeomorphisms.  \enddemo\qed

 \medbreak

We conclude this section by summarizing our principal conclusions.

\medbreak

\proclaim {Theorem 4.6 }  
Suppose that $A$ and $B$ are separable 
$C^*$-algebras and that $A$ is   $K$-nuclear.
 Then   there are natural
isomorphisms of topological groups
$$
\KKgraded 1.A.B  _S  \,\cong\,
\KKgraded 1.A.B  _Z  \,\cong\,
\KKgraded 1.A.B  _C  \,\cong\,
\KKgraded 1.A.B  _q    .
$$
\medbreak\flushpar
Thus  if all $C^*$-algebras appearing in the 
first variable of $KK$ are $K$-nuclear, then 
the Kasparov $KK$-pairing is
separately continuous in the Salinas topology.
\endproclaim\qed

{\bf {Remark 4.7}}  In our applications we normally assume that
$A$ is $K$-nuclear since our applications deal with Salinas's work
on quasidiagonality and hence require that the Kasparov pairing be
continuous in the topology that he uses. So {\bf{ we shall take
all $C^*$-algebras appearing in the first variable of $KK$ to be
$K$-nuclear for the rest of the paper.}} Generally speaking the
results still hold without this assumption, but they have to be
interpreted as holding for the  topology $\KKgraded *.A.B  _q $.

 \newpage\beginsection {5. Consequences of the
 continuity of the $KK$-product }

In this section we begin our application of the continuity of the
$KK$-pairing in order to show that other structural maps are
continuous. The first result is an immediate and very useful
consequence.

\medbreak

\proclaim {Theorem 5.1 }  Suppose that  the $C^*$-algebra
$D$ is $KK$-equivalent to $D'$
via an invertible $KK$-class $y \in \KKgraded 0.D.{D'}      $. Then for any
$C^*$-algebra $B$
the induced isomorphism
$$
y \tensor _{D'} ( - ) : \KKgraded *.{D'}.B \longrightarrow
 \KKgraded *.{D }.B
$$
is an isomorphism of topological groups
  and for any $C^*$-algebra $A$  the induced isomorphism
$$
(-)\tensor _D y  :  \KKgraded *.A.D   \to  \KKgraded *.A.{D'}
$$
is an isomorphism of topological groups.
\endproclaim

\medbreak

\demo {Proof}   Let $y' \in \KKgraded 0.{D'}.D $ be the $KK$-inverse
of $y$. Then the map
$$
y' \tensor _{D} ( - ) : \KKgraded *.{D}.B \longrightarrow
 \KKgraded *.{D' }.B
$$
is the inverse to the map $y \tensor _{D'} ( - )$. Both of these
maps are continuous, by Theorem 4.6, and that implies that the map
$y \tensor _{D'} ( - )$ is an isomorphism of topological groups.
The second part of the theorem is proved similarly.
 \qed
\enddemo
\medbreak

Suppose given a short exact sequence
$$
0  \to J \overset{i}\to\longrightarrow A \overset{p}\to\longrightarrow  A/J  \to
0
$$
which is  split  with splitting $s: A/J \to A$. There is a natural
$KK$-class
$$
[1_A] - [sp]  \in \KKgraded 0.A.A
$$
 and this class lies in the image of the injection
$$
\KKgraded 0.A.J  \overset{i_*}\to\longrightarrow \KKgraded 0.A.A
.
$$
We denote by $t^s \in \KKgraded 0.A.J $ the unique class
satisfying
$$
i_* (t^s) =  [1_A] - [sp ] .
$$

\medbreak

\proclaim {Proposition 5.2}
\roster
\item   Suppose that
$$
0  \to J \overset{i}\to\longrightarrow A \overset{p}\to\longrightarrow  A/J  \to
0
$$
is a split short exact sequence with splitting $s: A/J \to A$.
Then there is a canonical\footnote{determined uniquely by the $KK$-class of the
splitting $s: A/J \to A$}
 isomorphism of topological groups
$$
\KKgraded *.A.B \,\cong\, \KKgraded *.J.B \oplus \KKgraded *.{A/J}.B   .
$$
The continuous structural maps are given by
$$
i^*:   \KKgraded *.A.B        \longrightarrow        \KKgraded *.J.B   ,
$$
$$
s^*:    \KKgraded *.A.B     \longrightarrow     \KKgraded *.{A/J}.B    ,
$$
$$
p^* :   \KKgraded *.{A/J}.B         \longrightarrow       \KKgraded *.A.B    ,
$$
and
$$
t^s :  \KKgraded *.J.B  \longrightarrow  \KKgraded *.A.B   .
$$
\item   Suppose that
$$
0  \to J \overset{i}\to\longrightarrow B \overset{p}\to\longrightarrow  B/J  \to
0
$$
is a split short exact sequence with splitting $s: B/J \to B$.
Then there is a canonical isomorphism of topological groups
$$
\KKgraded *.A.B \,\cong\, \KKgraded *.A.J  \oplus \KKgraded *.A.{B/J}   .
$$
with continuous structural  maps which are analogous to  the maps of Part (1).
\endroster
\endproclaim
\medbreak

\demo {Proof} (1):  The algebraic isomorphism follows immediately
from exactness properties of the Kasparov groups, and three of the
four structural maps are obviously continuous. The real point is
that the map $t^s $  also is continuous, which follows from
Theorem 4.6.  The proof of (2) is similar. \qed\enddemo

Note that whenever there is an isomorphism $G \cong G_1\oplus G_2 $
of topological groups then the   projection maps are
open.
We use this fact in the following theorem.

\proclaim {Theorem 5.3} Suppose that $\{A_j\} $ is a countable
family of   $C^*$-algebras. Then the natural isomorphism
of   groups
$$
\varphi : \KKgraded *.{\oplus A_j}.B   \longrightarrow \prod \KKgraded *.{A_j}.B
$$
is an isomorphism of topological groups, where the right hand side is topologized
as the
  product of the topological groups $\KKgraded *.{A_j}.B  $.
\endproclaim

\medbreak
\demo {Proof} The map $\varphi $ is the product of the maps
$$
\varphi _k : \KKgraded *.{\oplus A_j}.B  \longrightarrow \KKgraded *.{  A_k}.B
$$
induced by the canonical inclusions $A_k \to \oplus A_j $, and the
map  $\varphi $ is known to be an isomorphism, by work of J.
Rosenberg  \cite {18, 1.12}.  Each map $\varphi _k $ is
continuous, and this implies that the map $\varphi $ is
continuous.

It remains to show that the map $\varphi $ is open, and for this it suffices
to show that each $\varphi _k$ is open.
 This is indeed the case, by Proposition 5.2,
since each sequence
$$
0 \to A_k \to \oplus A_j \to  (\oplus A_j)/A_k  \to  0
$$
is split exact.
\qed\enddemo

 \medbreak

\newpage\beginsection{6. Polonais and Pseudo\-polonais groups}

Recall from C.C. Moore \cite {16, 17}
 that a topological space $X$ is said to be {\it {polonais }} if it is complete, separable, and metric. If $X$ is a topological group then we also
insist that the metric be invariant under the group action.

The motivating example is the space $Hom(A,B)$, where $A$ and $B$
are separable. Then this is a polonais space.

\medbreak

\proclaim {Definition 6.1 } A topological group $G$ is said to be
 {\it {pseudo\-polonais}} if it is separable, pseudo\-metric, and
if the quotient group $\overline G = G/G_o$ of the group $G$ by
$G_o$, the closure of zero,  with the quotient topology
 is a polonais group.\endproclaim

We note that this is in the direction but not as general as the
weakening of the
 polonais hypothesis considered by C. C. Moore \cite {17, p. 10}.
\medbreak

If $A$ and $B$ are separable then the space $[A,B]$ with its
natural topology is separable, complete,
 and has a pseudo\-metric $\hat \mu $
given as in Section 4.

\proclaim{Theorem 6.2} The groups $\KKgraded *.A.B $ have
a natural structure as pseudo\-polonais topological groups
given by the isomorphism of topological groups
$$
\KKgraded 0.A.B  \cong [qA, B\tensor\Cal K ]
$$
$$
\KKgraded 1.A.B  \cong [q(SA), B\tensor\Cal K ] .
$$
\endproclaim

\medbreak
\demo{Proof} Theorem 4.6 gives us the isomorphisms as indicated. 
So it suffices to show that $[qA, B\tensor\Cal K] $ has the 
structure of a pseudopolonais topological group. Given the
observation above, the group is complete, separable, and 
pseudopolonais. It remains to check that the group action is
invariant under the metric. The group addition operation
is given as follows: if $f, g : qA \to B\tensor\Cal K $ then
$[f] + [g]$ is the homotopy class of the composite
$$
qA \overset{(f,g)}\to\longrightarrow 
(B\tensor \Cal K) \oplus (B\tensor\Cal K)\,\, \subseteq \,\, 
B\tensor\Cal K \tensor M_2(\Bbb C) \,\,\cong\,\, B\tensor\Cal K
$$
and it is easy to check that this respects the pseudometric.
\enddemo\qed

Next we recall some important facts about polonais groups which
are all found in \cite {17,  2.3}. \medbreak

\proclaim {Theorem 6.3 } \roster
\item If $\{G_i \}$ is a sequence of polonais groups then
their product $\prod _i G_i $ with the product topology is
polonais.
\medbreak
\item Suppose that
$$ 0 \to G' \overset{i}\to\longrightarrow  G
\overset{j}\to\longrightarrow G'' \to 0 $$ is a short exact
sequence of Hausdorff topological groups with $i$ a homeomorphism
onto its image and $j$ continuous and open. Then $G$ is polonais
if and only if $G'$ and $G''$ are polonais.
\item Let $G_1$ and $G_2$ be separable metric groups with $G_1$
polonais, and let $\phi : G_1 \to G_2$ be a Borel homomorphism.
Then $\phi $ is continuous.
\item Let $G_1$ and $G_2$ be polonais and let $\phi : G_1 \to G_2$
be a continuous bijection. Then $\phi $ is open and hence an
isomorphism of topological groups.
\endroster
\endproclaim

\medbreak

\demo {Proof} Part 1 is \cite {17, Prop. 2}
 and Part 2 is \cite {10, Prop. 3}.
Parts 3 and 4 appear as \cite {17, Prop. 5}; Moore  attributes
them to Banach \cite {1} and to Kuratowski \cite {15}
respectively. \qed\enddemo

\medbreak

Here is an easy consequence.

\medbreak \proclaim {Proposition 6.4  } Suppose that  $G$ and $H$
are polonais groups and  that
 $f: G \to H$ is a continuous homomorphism.
\roster
\item
Suppose that   $Im(f)$ is a closed subgroup of $H$. Give
$G/Ker(f)$ the quotient topology from G and $Im(f)$ the relative
topology from $H$. Then the natural bijection $$ \hat {f}  :
G/Ker(f) \to Im(f) $$ is an isomorphism of topological groups
and hence $f$ itself is relatively open.
\medbreak
\item If $f$ is onto then $\hat f : G/Ker(f) \to H $
is an isomorphism of topological groups.
\item If $$ G_1 \overset{f_1}\to\longrightarrow  G_2
 \overset{f_2}\to\longrightarrow G_3
$$ is an exact sequence of polonais groups then $Im(f_1) \cong
Ker(f_2) $ is closed, and hence,
 topologizing as in 1), the natural map is an isomorphism
$$ G_1/Ker(f_1) \cong  Im(f_1) $$ of topological groups.
\medbreak

\endroster
\endproclaim
\medbreak \demo{Proof} Both $G/Ker(f)$ and $Im(f)$ are polonais in
their specified topologies, and $\hat f$ is a continuous
bijection. So 1) is  immediate from Theorem 6.3 (4). Parts 2) and
3) follow at once from 1). \qed \enddemo

\medbreak \medbreak

\proclaim {Proposition 6.5} Suppose that $G$ and $H$ are
pseudo\-polonais groups and that $\beta : G \to H$ is a continuous
homomorphism inducing $\beta ' : G_o \to H_o$. Suppose that
$\overline \beta: \overline G \to \overline H $ is an open map.
 \roster
\item If both $\beta $
and $\beta '  $ are  onto,
  then $\beta $ is an open map.
\medbreak
\item If both $\beta $
and $\beta '  $ are  bijections,
  then $\beta $ is an isomorphism of topological groups.
\endroster
 \endproclaim

 \medbreak

 \demo {Proof}

We first consider 1). Let $U$ be an open neighborhood of $0$ in
$G$. Then of course $G_0 \subseteq U$. We wish to show that $\beta
(U)$ is open in $H$. Note that 
$$
\beta ^{-1}\beta U = \{ u + z : u \in U, \,\, z \in Z \}
= \bigcup _{z\in Z} (U + z)
$$
and this is a union of open sets, hence  open.
Thus
without loss of generality, we may assume
that $U$ is 
open and
saturated   with respect to $\beta $;   i.e.,  $U = \beta
^{-1}\beta U$. Consider the commuting diagram $$ \CD G   @>\beta
>>   H   \\ @VV\pi V     @VV\pi V  \\ \overline G    @>{\overline\beta
}>>   \overline H.
\endCD
$$ The map $\overline\beta $ is open by assumption, and hence
$\pi\beta U$ is open in $\overline H$. This implies that the $\pi
$-saturation of $\beta U$, is open in $H$. However, $$ \align
Sat(\beta U) &= \beta (U) + H_0  \\ &= \beta (U) + \beta (G_o)
\\ &=\beta ( U +  G_o)    \\ &= \beta (U)
\endalign
$$ and this completes the proof of 1). Part 2) follows from 1).
 \qed\enddemo

\medbreak

\proclaim {Corollary 6.6} Suppose that $G$ is a pseudo\-polonais
group with closed subgroup $G'$. Let $G''$ be a polonais  group
and suppose that there is a short exact sequence $$ 0 \to G'  \to
G  \overset{f}\to\longrightarrow G'' \to 0 $$ with $f$ continuous.
Give $G/G'$ the quotient polonais topology. Then the natural
algebraic isomorphism $$ \hat f:  G/G'
\overset\cong\to\longrightarrow G'' $$ is an isomorphism of
topological groups.
\endproclaim

\medbreak \demo {Proof} This is immediate from 6.3 and the fact
that polonais groups are Hausdorff. \qed
\enddemo

\medbreak

Our proof of the following result from folklore  is based entirely
upon ideas of C. C. Moore, particularly \cite {17, Prop. 11} and
\cite {16, Prop. 1.4}.

\medbreak

\proclaim {Theorem 6.7} Suppose that $G$, $H$, and $K$ are
polonais groups and that
$$
m : G \times H \longrightarrow K
$$
is a bilinear pairing which is separately continuous.
Then $m$ is jointly continuous.
\endproclaim
\medbreak

\demo {Proof} The fact that $m$ is separately continuous implies
that the map $m$ is jointly measurable, by a result of Kuratowski
\cite {15, p. 285}. A theorem of Banach \cite {1, p.  25} implies
that there exists a set $P \subseteq G \times H $ of first
category such that the function
$$
m : (G \times H)  -  P   \longrightarrow K
$$
is continuous.

     We write $m(g,h) = gh$ for convenience.
Suppose that $(g_n, h_n ) \to (g_o, h_o) $ in $G \times H $.
We must prove that $g_nh_n  \to g_oh_o$ in $K$.
Consider the set
$$
\tilde P = \cup _{n=1}^\infty  P(-g_n,-h_n)  .
$$
This is a set of first category in $G \times H$, and so
$\tilde P \neq G \times H$.
Choose  some $(g,h) \in (G \times H) - \tilde P$.
Then
$$
(g + g_n,h +  h_n )  \in   (G\times H)  - P  \qquad \forall n
$$
and since $m$ is continuous on that set, we have
$$
(g + g_n)(h +  h_n ) \to (g + g_o)(h +  h_o)  .
$$
We know that $g_nh \to g_oh$ and that $gh_n \to gh_o$
and so subtracting these terms and the constant term
yields
$
g_nh_n \to g_oh_o $
as required.\qed \enddemo

\medbreak

At last we can complete the proof of the main theorem of this paper.

\medbreak

\proclaim{Theorem 6.8} The Kasparov pairing
$$
\KKgraded *.{A_1}.{B_1\tensor D} \times \KKgraded *.{D\tensor A_2}.{B_2}
\overset{\tensor _D}\to\longrightarrow
\KKgraded *.{A_1\tensor A_2}.{B_1\tensor B_2}
$$
is jointly continuous in the Salinas topology, provided that
all $C^*$-algebras appearing in the first variable are $K$-nuclear.
\endproclaim
\medbreak

\demo{Proof} For simplicity, let
$$
\align
G &= \KKgraded *.{A_1}.{B_1\tensor D}  \\
H &=  \KKgraded *.{D\tensor A_2}.{B_2}  \\
K &=  \KKgraded *.{A_1\tensor A_2}.{B_1\tensor B_2} \\
m(x,y) &= x\tensor _D y
\endalign
$$
Then we must show that $m: G \times H \to K$ is jointly continuous.
We know that $m$ is separately continuous, by Theorem 4.6.
Recall that $G_o$ denotes
 the closure of zero in $G$,   $\overline G = G/G_o$,
and similarly for $H$ and $K$. The map $m$ induces a natural
map
$$
\overline m : \overline G \times \overline H \to \overline K
$$
and it is easy to see that this map is also separately continuous.
Now the groups $\overline G$, $\overline H$, and $\overline K$
are polonais, by Theorem 6.2. Thus we may apply Theorem 6.7
to conclude that the map $\overline m$ is jointly continuous.

Let $\pi $ generically denote the map from a group to its Hausdorff
 quotient. Since the diagram
 $$
 \CD
 G \times H  @>{m}>>   K   \\
 @VV{\pi \times \pi}V   @VV{\pi}V  \\
 \overline G \times \overline H  @>{\overline m}>>  \overline K
 \endCD
 $$
commutes, the
 composite
 $$
 G \times H \overset{m}\to\longrightarrow K
 \overset\pi\to\longrightarrow \overline K
 $$
 is jointly continuous. 
 Finally, let $U$ be an open neighborhood of $0 \in K$.
   The map 
   $$
   \pi : K \to \overline K
   $$
    is
  open for any pseudo\-polonais group, and hence $\pi U$ is
  open in $\overline K$. Then
  $$
  m^{-1}(U) = (\pi m)^{-1}(\pi U) =
  (\pi \times \pi)^{-1}\overline m^{-1}(\pi U)
  $$
  which is open since both $\overline m$ and $\pi \times \pi $ are
  continuous.
  Thus $m$ is jointly continuous and so the proof is complete.
  \qed\enddemo

\medbreak

{\bf {Remark 6.9}}  \,\, The study of the topological structure of
$\KKgraded *.A.B $ was initiated by L. G. Brown, R.G. Douglas, and
P.A. Fillmore
  in their 1973    \cite {5} announcement.\footnote
{At that time it was not known that
$$ \Cal Ext _*(A) \cong \KKgraded *.A.{\Bbb C} \equiv
K^*(A),
$$
 but I shall use the $K^*$ notation for the
convenience of the reader.}  They considered the natural topology
on $K^*(A)$ (usually they concentrated upon $K^*(C(X))$; but   it
was understood   that commutativity was inessential), noted that
$0$ was not necessarily closed in this topology, and used the fact
that $0$ {\it{was}} closed for $X \subset  \Bbb C$ to demonstrate
that the set of bounded operators of the form {\it {(normal) +
(compact)}} was norm-closed.

They also introduced the group $PExt(X)$, defined as
$$
PExt(X) = \{ x \in K^1(C(X)) : f_*x = 0 \,\,\, \forall\, f \}
$$
 where $f$ ranges
over all continuous functions from $X$ to an ANR.\footnote {The
\lq\lq P\rq\rq\,\, stood for \lq\lq pathological\rq\rq, not \lq\lq
pure\rq\rq , and thus the notation is serendipitous.} Writing $X =
\invlim X_n$ where the $X_n$ are finite complexes, they announced
that
$$
 PExt(X) = Ker \,\, [ K^1(C(X))
\to \invlim K^1(C(X_n)) ].
$$
 Given Bott periodicity, the
Milnor $\invlimone $ sequence, and the UCT (none of these were
completely established in 1973),  this was tantamount to proving
that
 $$ PExt(X) \,\,\cong  \,\, \pext {K^0(X)}.{\Bbb Z}   .
 $$
\medbreak\flushpar Finally, they announced that
$$
 Ker\,\,[\, \gamma _{\infty }  :
K^1(C(X)) \to \hom {K^1(X)}.{\Bbb Z} \,]
 $$
 was the maximal
compact subgroup of $K^1(C(X))$. All of these results were
announced in \cite {5} but were not discussed in subsequent
publication, except for the case $X \subset \Bbb C$.

It is particularly gratifying to finally establish that the
Kasparov groups are pseudo\-polonais. L.G. Brown asked about this
matter in the early days of the BDF groups, and he has asked the
author regularly every few years since. His persistence is
appreciated. We should note that his original question asked if it
was possible to write the groups as the quotient of one polonais
group by another. This question is still open. \medbreak

\newpage\beginsection {7. Suspensions, Images, Boundary Maps,
and the Index Map  }

 In this section we show that various
suspension and boundary maps are continuous.
We also  show that the boundary homomorphisms in $KK$-long exact sequences
are given by instances of the $KK$-pairing and hence
are continuous in each variable. Finally, we demonstrate that
the index map
$$
\gamma : \KKgraded *.A.B \longrightarrow \hom {K_*A)}.{K_*(B)}
$$
is continuous in the natural topology on $ Hom $.

\medbreak

The following proposition is well-known and we include a proof
only for convenience.
\medbreak

\proclaim {Proposition 7.1} Suppose given a short
exact sequence of $C^*$-algebras
$$
0 \to J \overset{i}\to\longrightarrow C
 \overset{p}\to\longrightarrow C/J \to 0   .
$$
Then:
\roster
\item If  $C/J$ contractible, then  $[i] \in \KKgraded 0.J.C    $ is
  $KK$-invertible.
\medbreak
\item If $J$ is contractible, then $[p]  \in \KKgraded 0.C.{C/J}   $ is
  $KK$-invertible.
\endroster
\endproclaim
\medbreak
\demo {Proof}
 Suppose that $C/J $ is contractible.
Then the map
$$
i_* : \KKgraded 0.C.J \to \KKgraded 0.C.C
$$
is an isomorphism,
so define $\sigma \in \KKgraded 0.C.J $ to be the unique element
with
$$
\sigma \tensor _J [i] \,=\,  i_*\sigma = [1_C]  .
$$
 Consider the commutative diagram
$$
\CD
\KKgraded 0.C.J   @>{i^*}>>   \KKgraded 0.J.J   \\
@VV{i_*}V    @VV{i_*}V    \\
\KKgraded 0.C.C   @>{i^*}>>   \KKgraded 0.J.C
\endCD
$$
Computing, we have
$$
i_*i^*\sigma = i^*i_*\sigma = i^*[1_C] =
[i] = i_*[1_J]
$$
and since $i_*$ is an isomorphism, it follows that
$$
[i] \tensor _C \sigma \,=\, i^*\sigma \,=\, [1_J]
$$
so that
$\sigma \in \KKgraded 0.C.J $ is a $KK$-inverse to $[i] \in \KKgraded 0.J.C $.
This demonstrates Part (1).

The proof of Part (2) is similar and uses $\sigma \in \KKgraded 0.{C/J}.C $
with $p^*\sigma = [1_C] $. \qed\enddemo

\medbreak

For any separable $C^*$-algebra $E$, let
$$
b_E \in \KKgraded 1.E.{SE}
$$
be the Bott periodicity element. This element is the image of
the universal Bott element
$$
b \in \KKgraded 1.{\Bbb C}.{S\Bbb C} \,\,\cong\,\, \KKgraded 0.{\Bbb C}.{\Bbb C}
\,\cong\, \Bbb Z
$$
under the canonical structural map
$$
 \KKgraded 1.{\Bbb C}.{S\Bbb C} \to
\KKgraded  1.{\Bbb C\tensor E}.{S\Bbb C\tensor E} \cong \KKgraded 1.E.{SE}   .
$$

\medbreak

Much of the following proposition is due to Kasparov \cite {14, p.
566}.

\medbreak
\proclaim {Proposition 7.2} Suppose that
$$
0 \to D' \to D \overset{f}\to\longrightarrow  D'' \to 0
\tag *
$$
is a short exact sequence of $K$-nuclear $C^*$-algebras. Then there exists
a canonical class
$$
\Delta \in \KKgraded 1.{D''}.{D'}
$$
such that for all $C^*$-algebras $A$ the boundary homomorphism
$$
\,\delta _* : \KKgraded *.A.{D'' }  \to \KKgraded {*-1}.A.{D'}
$$
is given by
$$
\delta _*(x) = x \tensor _{D''}\Delta
$$
and for all $C^*$-algebras $B$ the boundary homomorphism
$$
\,\delta ^*: \KKgraded *.{D'}.B   \to \KKgraded {*-1}.{D'' }.B
$$
is given by
$$
\delta ^* (y) = \Delta \tensor _{D'} y .
$$
Both boundary maps
are continuous.
     If (*) is an essential extension then the class $\Delta $ corresponds to the
class of
the extension
(*) under the
identification
$$
 \Cal Ext(D'',D')  \overset\cong\to\longrightarrow \KKgraded 1.{D''}.{D'}    .
$$
Finally, if $D$ is contractible then the class $\Delta $ is $KK$-invertible.
\endproclaim
\medbreak

\demo {Proof}  The map $f: D \to D''$  has mapping cone sequence
$$
0 \to   SD''     \overset\zeta\to\longrightarrow
     Cf     \to     D     \to 0
$$
 and there is an associated  homotopy-commutative
diagram  with exact rows and columns
 \medbreak
$$
\CD
@.  @.  0   @.   @.   \\
@.   @.   @VVV   @.   @.  \\
@.   @.   D'    @.   @.   \\
@.   @.   @VV{\eta }V   @.   @.   \\
0     @>>>      SD''       @>{\zeta }>>       Cf          @>>>   D @>>>       0
     \\
@.   @VV{id}V            @VVV      @VV{f}V      @.            \\
0  @>>>    SD''     @>>>    CD''      @>>>     D''     @>>>      0   \\
@.   @.   @VVV   @.   @.  \\
@.  @.  0   @.   @.
\endCD
$$
 \medbreak
\flushpar  by \cite {20, 2.3}, where $CD'' $, the cone of $D'' $,
is contractible, and hence the natural inclusion $\eta : D' \to Cf
$ is  $KK$-invertible, by Proposition 5.1. Let $\eta ^{-1} \in
\KKgraded 0.Cf.D' $ denote its $KK$-inverse.
 There is then a commuting diagram
$$
\CD
\KKgraded {* }.A.{D'' }    @>{\,\delta _*}>>    \KKgraded {*-1}. A.{D'}  \\
@VV{(-)\tensor _{D''}b_{D''} }V      @AA{  (-)\tensor _{Cf}(\eta  ^{-1})}A    \\
\KKgraded {*-1}.A.{SD'' }    @>{(-)\tensor _{SD''}[\zeta ]}>>    \KKgraded
{*-1}.A.{Cf}
\endCD
$$
\medbreak
\flushpar and each vertical map is an isomorphism of topological groups.
Define
$$
\Delta = b_{D''}  \tensor _{SD''}   [\zeta ] \tensor _{Cf} (\eta ^{-1}) \in
\KKgraded 1.D'.D''  .
$$
Note that $\Delta $ is determined uniquely by the class $[f] \in
\KKgraded 0.D.{D''} $. The commutativity of the diagram implies
that $\delta _* (x) = x \tensor _{D''} \Delta $ as required. Since
$KK$-pairing  with any element is continuous, by Theorem 6.8, it
follows that the boundary homomorphism   is continuous. This
establishes the proposition in the second variable.

The argument in the first variable is similar with an additional twist.
The map $\delta ^*$ may be represented as the composite

$$
\CD
\KKgraded *.{D'}.B   @>{\delta ^*}>>  \KKgraded {* - 1}.{D''}.B  \\
@VV{(\eta ^{-1})\tensor _{D'}}(\,\,\,)V   @AA{b_{D''} \tensor _{SD''}(\,\,\,)}A
\\
\KKgraded *.{Cf}.B   @>{ [\zeta ] \tensor _{Cf} (\,\,\,) }>>  \KKgraded
*.{SD''}.B
\endCD
$$
\medbreak\medbreak

\flushpar so that
$$
\delta ^*(y) =  \Delta \tensor _{D'} y
$$
and so $\delta ^* $ is continuous as claimed.

Finally, suppose that $D$ is contractible. Then $\zeta $ is
$KK$-invertible, by Proposition 7.1(1). Since the Bott element and
$\eta ^{-1} $ are also $KK$-invertible and the $KK$-pairing is
associative, this implies that $\Delta $ is $KK$-invertible.
\qed\enddemo

\medbreak

{\bf {Remark 7.3}}  We note as a general statement
 that all $KK$-homology and cohomology operations
are given by appropriate $KK$-products and hence are continuous. For
instance, fix a positive integer $n$ and let $n : S\Bbb C \to S\Bbb C $
be the standard map of degree $n$. Then the mapping cone sequence has the
form
$$
0 \to S^2\Bbb C \to Cn  \to S\Bbb C \to 0
$$
and applying the functor $K_*(A\tensor (-))$ to this sequence
 gives rise in a manner that we have described elsewhere \cite {15}
to the Bockstein operation
$$
\beta _n : K_j(A; \Bbb Z/n ) \to K_{j-1}(A)  .
$$
Now this map is of course continuous, since the groups themselves are discrete.
A more interesting operation arises by applying
the functor $\KKgraded *.A.{B\tensor (-) } $.
Define
$$
KK_*(A,B; \Bbb Z/n)  = { \KKgraded *.A.{B\tensor Cn}  } .
$$
Then there is a Bockstein operation
$$
\beta _n : KK_j(A,B; \Bbb Z/n)  \to  \KKgraded {j-1}.A.B
$$
and it is a non-trivial fact that this map is also continuous. Its
continuity follows directly from 6.8.

\medbreak

The Kasparov pairing gives a natural index map
$$
\gamma : \KKgraded *.A.B  \to \hom {K_*(A)}.{K_*(B)}
$$
defined by
$$
\gamma (y) (x) = x \tensor _A y
$$
where the natural identifications
$$
\KKgraded *.{\Bbb C}.A  \cong  K_*(A)  \qquad
\KKgraded *.{\Bbb C}.B  \cong  K_*(B)
$$
are made without further
comment. We wish to topologize $ \hom {K_*(A)}.{K_*(B)} $
so as to be consistent with the Salinas topology,
   which is based upon the point-norm topology for extensions,
 and  so the natural
way to do this is to
regard $K_*(A)$ and $K_*(B)$ as discrete
 and to
use the topology of pointwise convergence on  $ \hom {K_*(A)}.{K_*(B)} $.
Assume this topology as given henceforth. Note that under this topology
the group
$$
 \hom {K_*(A)}.{K_*(B)}
$$
 is polonais: it is a closed subset
of the polonais, totally disconnected  group
$$
\prod _1^\infty K_*(B)
$$
 where $K_*(B)$
is given the discrete topology and the product topology is used
on the product. If $K_*(A)$ is finitely generated then
 $ \hom {K_*(A)}.{K_*(B)} $ is discrete, but in general this is not
the case.  For instance,
$$
\hom {\oplus _1^\infty  (\Bbb Z/2)}.{\,\Bbb Z/2} \,\,\cong\,\,
 \prod _1^\infty \,\hom  {\Bbb Z/2}.{\Bbb Z/2}  \,\,\cong\,\,
 \prod _1^\infty \Bbb Z/2
$$
is  topologically a Cantor set. In general, $\hom G.H $ is
totally disconnected (I am endebted to George Elliott for this
point).

\medbreak
\medbreak

\proclaim {Proposition 7.4} The natural map
 $$
\gamma : \KKgraded *.A.B  \to \hom {K_*(A)}.{K_*(B)}
$$
is continuous.
If $Im(\gamma )$ is closed then $\gamma $ is open onto its
image.
  If $\gamma $ is an algebraic isomorphism
then $\gamma $ is an isomorphism of topological groups.
\endproclaim
\medbreak

\demo {Proof} Suppose that $y^\alpha $ is a net
in $\KKgraded *.A.B $ which converges to $y \in \KKgraded *.A.B $.
Then for each $x \in K_*(A) $,
 $$
\gamma (y^\alpha )(x) =   x\tensor _A y^\alpha
$$
 is a net
in $K_*(B) $ which converges to
$$
\gamma (y)x = x\tensor _B y \in K_*(B)
$$
since the map $x \tensor _B (-) $ is continuous by Theorem 4.6.
Thus $\gamma (y^\alpha )$ converges pointwise to $\gamma (y) $ as
desired. So $\gamma $ is continuous.

Now suppose that $Im(\gamma ) $ is closed,
so that it is polonais. It suffices to prove that $\gamma $ is open.
 Factor the map $\gamma $
as $$ \KKgraded *.A.B  \overset\pi\to\rightarrow \bKKgraded *.A.B
\overset{\overline\gamma}\to\rightarrow Im(\gamma )   . $$ Since
$\pi $  is an open map, it suffices to show that the map
$\overline\gamma $ is open. However the group $\bKKgraded *.A.B $
is polonais, and $Ker(\overline\gamma )$ is a closed subgroup, so
the quotient group $$ \KKgraded *.A.B /Ker(\gamma ) = \bKKgraded
*.A.B /Ker(\overline\gamma ) $$ is also polonais. Thus the induced
map $$ \hat\gamma :  \KKgraded *.A.B /Ker(\gamma ) \to Im(\gamma )
$$ is a continuous bijection of polonais groups, and hence a
homeomorphism, by Theorem 6.3(4), and of course this implies that
$\overline\gamma $ is open. \qed\enddemo

\medbreak

\end